\documentclass[3p,,preprint,12pt]{elsarticle}
\makeatletter\if@twocolumn\PassOptionsToPackage{switch}{lineno}\else\fi\makeatother

\usepackage{tabulary,xcolor}
\usepackage{amsfonts,amsmath,amssymb}
\usepackage[T1]{fontenc}
\usepackage{tikz}
\usepackage{flowchart} \usetikzlibrary{arrows}
\usetikzlibrary{shapes, arrows.meta, positioning}
\makeatletter
\let\save@ps@pprintTitle\ps@pprintTitle
\def\ps@pprintTitle{\save@ps@pprintTitle\gdef\@oddfoot{\footnotesize\itshape \null\hfill\today}}
\def\hlinewd#1{%
  \noalign{\ifnum0=`}\fi\hrule \@height #1%
  \futurelet\reserved@a\@xhline}

\AtBeginDocument{\ifNAT@numbers \biboptions{sort&compress}\fi}
\makeatother

  \renewenvironment{abstract}{\global\setbox\absbox=\vbox\bgroup
    \hsize=\textwidth%
  \noindent\unskip\textbf{}
   \par\medskip\noindent\unskip\ignorespaces}
   {\egroup}
  
\usepackage{ifluatex}
\ifluatex
\usepackage{fontspec}
\defaultfontfeatures{Ligatures=TeX}
\usepackage[]{unicode-math}
\unimathsetup{math-style=TeX}
\else 
\usepackage[utf8]{inputenc}
\fi 
\ifluatex\else\usepackage{stmaryrd}\fi

\usepackage{url,multirow,morefloats,floatflt,cancel,tfrupee}
\makeatletter

\AtBeginDocument{\@ifpackageloaded{textcomp}{}{\usepackage{textcomp}}}
\makeatother
\usepackage{colortbl}
\usepackage{xcolor}
\usepackage{pifont}
\usepackage[nointegrals]{wasysym}
\makeatletter
\usepackage{graphicx}
\usepackage{amsmath}      
\usepackage{capt-of}
\usepackage{mathptmx}
\usepackage{float}
\usepackage{listings}
\def\mcWidth#1{\csname TY@F#1\endcsname+\tabcolsep}

\def\cAlignHack{\rightskip\@flushglue\leftskip\@flushglue\parindent\z@\parfillskip\z@skip}
\def\rAlignHack{\rightskip\z@skip\leftskip\@flushglue \parindent\z@\parfillskip\z@skip}

\usepackage{ifxetex}
\ifxetex\else\if@twocolumn\usepackage{dblfloatfix}\fi\fi

\AtBeginDocument{
\expandafter\ifx\csname eqalign\endcsname\relax
\def\eqalign#1{\null\vcenter{\def\\{\cr}\openup\jot\m@th
  \ialign{\strut$\displaystyle{##}$\hfil&$\displaystyle{{}##}$\hfil
      \crcr#1\crcr}}\,}
\fi
}

\AtBeginDocument{%
  \@ifpackageloaded{endfloat}%
   {\renewcommand\efloat@iwrite[1]{\immediate\expandafter\protected@write\csname efloat@post#1\endcsname{}}}{\newif\ifefloat@tables}%
}%

\def\BreakURLText#1{\@tfor\brk@tempa:=#1\do{\brk@tempa\hskip0pt}}
\let\lt=<
\let\gt=>
\def\processVert{\ifmmode|\else\textbar\fi}

\@ifundefined{subparagraph}{
\def\subparagraph{\@startsection{paragraph}{5}{2\parindent}{0ex plus 0.1ex minus 0.1ex}%
{0ex}{\normalfont\small\itshape}}%
}{}

\newcommand\role[1]{\unskip}
\newcommand\aucollab[1]{\unskip}
  
\@ifundefined{tsGraphicsScaleX}{\gdef\tsGraphicsScaleX{1}}{}
\@ifundefined{tsGraphicsScaleY}{\gdef\tsGraphicsScaleY{.9}}{}
\def\checkGraphicsWidth{\ifdim\Gin@nat@width>\linewidth
	\tsGraphicsScaleX\linewidth\else\Gin@nat@width\fi}

\def\checkGraphicsHeight{\ifdim\Gin@nat@height>.9\textheight
	\tsGraphicsScaleY\textheight\else\Gin@nat@height\fi}

\def\fixFloatSize#1{}
\let\ts@includegraphics\includegraphics

\def\inlinegraphic[#1]#2{{\edef\@tempa{#1}\edef\baseline@shift{\ifx\@tempa\@empty0\else#1\fi}\edef\tempZ{\the\numexpr(\numexpr(\baseline@shift*\f@size/100))}\protect\raisebox{\tempZ pt}{\ts@includegraphics{#2}}}}

\AtBeginDocument{\def\includegraphics{\@ifnextchar[{\ts@includegraphics}{\ts@includegraphics[width=\checkGraphicsWidth,height=\checkGraphicsHeight,keepaspectratio]}}}

\DeclareMathAlphabet{\mathpzc}{OT1}{pzc}{m}{it}

\def\URL#1#2{\@ifundefined{href}{#2}{\href{#1}{#2}}}

\def\UrlOrds{\do\*\do\-\do\~\do\'\do\"\do\-}%
\g@addto@macro{\UrlBreaks}{\UrlOrds}

\edef\fntEncoding{\f@encoding}

\makeatother

\begin{document}
	
	\begin{frontmatter}
		
		\title{\mbox{}Fibonacci Neural Network Approach for Numerical Solutions of Fractional Order Differential Equations
		}
		
		\author{Kushal Dhar Dwivedi*$^{1}$}
		\ead{kddwivedi1993@gmail.com, kddwivedi@iiitk.ac.in}
		\author{Anup Singh$^{2}$}
		\ead{anup.singh254@gmail.com}
		\author{Anirban Majumdar$^{1}$}
		\ead{anirban@iiitk.ac.in}
		\address
		{${}^1$Department of Sciences IIITDM, Kurnool-518007, India} 
		
		\address{${}^2$Department of Mathematics,\\ Institute of Technology Nirma University, Ahmedabad-382481, Gujarat, India. }
	\begin{abstract}
		\textbf{Abstract }
		In this paper, the authors propose the utilization of Fibonacci Neural Networks (FNN) for solving arbitrary order differential equations. The FNN architecture comprises input, middle, and output layers, with various degrees of Fibonacci polynomials serving as activation functions in the middle layer.
		 The trial solution of the differential equation is treated as the output of the FNN, which involves adjustable parameters (weights).
		  These weights are iteratively updated during the training of the Fibonacci neural network using backpropagation. The efficacy of the proposed method is evaluated by solving five differential problems with known exact solutions, allowing for an assessment of its accuracy. Comparative analyses are conducted against previously established techniques, demonstrating superior accuracy and efficacy in solving the addressed problems.
	\end{abstract} 
	\end{frontmatter}

\section{Introduction}
Differential equations have long served as a foundational tool for describing diverse physical phenomena across engineering, mathematics, physics, and economics. Despite their utility, obtaining solutions for nonlinear differential equations, which often govern these phenomena, poses a significant challenge. The complexity inherent in many of these equations frequently precludes exact analytical solutions. Consequently, numerical methods have emerged as indispensable tools for tackling such equations. These numerical approaches offer avenues for approximating solutions to differential equations, enabling the exploration and understanding of a wide range of real-world problems. Many researchers have contributed to developing efficient and accurate numerical methods to solve differential equations viz, Runge Kutta, Finite difference, Finite element methods, and so on. In the modern era, researchers are focused not only on developing a numerical method, but also it performs better than previously existing methods. Zhao et al. \cite{zhao2006new} solved backward stochastic differential equations by discretizing them in time-space discrete grids and then used the Monte Carlo scheme to approximate mathematical expectations. Rehman et al. \cite{ur2012numerical} employed Haar wavelet operational matrices for solving differential equations very efficiently. Yan et al. \cite{yan2014higher} have solved the fractional order differential equation(FDEs) by direct discretization of the fractional differential operator and discretization of the integral form of the FDEs. Li et al. \cite{li2017new} solved variable-order differential equations by reproducing kernel method. Baleanu et al. \cite{baleanu2018chebyshev} solved fractional differential equations with Chebyshev operational matrix. etc. Shiri et al. \cite{shiri2019system} in 2019, solved the system of FDEs, and the same author \cite{shiri2020collocation} in 2020 solved tempered FDEs with collocation method. Khiabani et al. \cite{dadkhah2020spline} used spline collocation method to solve fractional models. Dadkhah et al. \cite{dadkhah2020visco} modeled the FDEs for Visco-elastic dampers and solved it with spline collocation methods.\par 
 Artificial neural networks (ANN) are gaining attention due to their capability of doing outstanding tasks, viz., fitting into complex data to predict future outcomes, speech recognition, building virtual assistants like Google Home and Alexa, deep fake, and image conversion. It even helps machines become more productive and intelligent. This is all possible because the building block of ANN is to mimic the human brain. This is why ANN performs more efficiently than humans in many areas with no errors. The present work is just a contribution to show that ANN can be used to solve fractional differential equations (FDEs) more efficiently than previously existing methods. Due to its popularity, many authors have used the ANN for solving the FDEs in the last few years, viz, Balasubramaniam et al. \cite{balasubramaniam2006solution} solved the matrix Riccati differential equation by using ANN. Reynaldi et al. \cite{reynaldi2012backpropagation} solved the differential equation and the inverse problem of differential equation with ANN. Jafarian et al. \cite{jafarian2017artificial} used ANN to solve a class of FDEs. Z{\'u}{\~n}iga-Aguilar et al. \cite{zuniga2017solving} solved variable-order differential equation by using ANN with Mittag-Leffler kernel. \par 
 
In this study, the authors introduce a novel approach developed using FNN, which incorporates various degrees of Fibonacci polynomials as activation functions in the hidden layer. The FNN architecture consists of one input layer, one hidden layer, and one output layer. Each perceptron in the hidden layer is constructed using different degrees of Fibonacci polynomials with unit weights. The output layer aggregates the outputs of all perceptrons after applying different weights. Initially, these weights are assigned randomly, and then they are updated using an appropriate backpropagation algorithm during the training process, which will be elaborated on in subsequent sections.
The motivation behind using the FNN lies in approximating the desired solution of the FDEs by assigning the output of FNN as the solution, which is constructed using different degrees of Fibonacci polynomials. After that, residual will be formed with the help of considered FDE and its initial or boundary conditions. Lastly, the aim is to minimize the residual to get the solution of the taken FDES.
The proposed method is applied to various differential equations, demonstrating its effectiveness compared to existing schemes through comprehensive comparisons and analyses.
\section{Preliminaries}

\subsection{Fractional order derivative in Caputo sense}
\label{sec:2}
The fractional order derivative of the function $f(t)$ in Caputo sense is defined as \cite{podlubny1998fractional}
\begin{equation*}
	({}^C_0D^\alpha_t f)(t)=\frac{1}{\varGamma(n-\alpha)} \int_{0}^{t}(t-\tau)^{\normalsize(n-\alpha-1)} f^{(n)}(\tau) d\tau,\hspace{0.2cm} \alpha>0,\tau>0,
\end{equation*}
where $n-1<\alpha<n,\hspace{0.2cm} n\in N$,\\
The Caputo fractional order derivative of order $\alpha$ of the polynomial $t^k$ is

\begin{equation}
	\label{eq1}
	{}^C_0D^\alpha_t t^k =
	\begin{cases}
		0\hspace{2.1cm}k\in{0,1,2,...,\lceil\alpha\rceil},\\
		\frac{\varGamma(k+1)}{\varGamma(k+1-\alpha)} t^{k-\alpha},\hspace{0.2cm}  k\in N,\hspace{0.2cm}k\geq \lceil\alpha\rceil,
	\end{cases}
\end{equation}	 	

where symbol $\lceil\cdot\rceil$ indicates the ceiling function, which is expressed by the set of equation $\lceil\alpha\rceil=\text{min}\{n\in Z|n\geq\alpha\},$ $Z$ represents the set of integers.

\subsection{Characteristic of Fibonacci Polynomial}
The Fibonacci polynomial of arbitrary order can be obtained with the help of following recurrence relation 
\begin{center}
	$F_{m+2}(x)=xF_{m+1}(x)+F_m(x)$,\hspace{0.3cm}$m\geq0$,
\end{center} 
with initial conditions (IC) as
\begin{center}
	$	F_0(x)=0,\hspace{0.3cm}F_1(x)=1.$
\end{center}
From the above relations we can obtain the Fibonacci polynomial of order $m$ as
\begin{equation}
	\label{eq2}
	F_m(x)=\sum_{r=0}^{\lfloor\frac{m-1}{2}\rfloor} {m-r-1\choose r} x^{m-2r-1} \hspace{0.2cm},
\end{equation}
where $\lfloor\cdot \rfloor$ indicates the floor function, which is expressed by the relation $\lfloor\alpha \rfloor=\text{max}\{n\in Z|n\leq\alpha\}$, $Z$ denotes the set of integer numbers.\\
Above discussed equation can be reformulated as
\begin{equation}
	\label{eq3}
	F_i(x)=\underset{(j+i)=odd}{\sum_{j=0}^{i}} \frac{(\frac{i+j-1}{2})!}{j! (\frac{i-j-1}{2})!} x^j,\hspace{0.4cm}i\geq0.
\end{equation}

Fractional order derivative of the Fibonacci polynomial of degree $i$ in Caputo sense can be obtain with the help of above equations \eqref{eq1} and \eqref{eq3} as

\begin{align}
	{}^C_0D^\alpha_x F_i(x) = &	{}^C_0D^\alpha_x \Big(\underset{(j+i)=odd}{\sum_{j=0}^{i}} \frac{(\frac{i+j-1}{2})!}{j! (\frac{i-j-1}{2})!} x^j \Big),\nonumber\\
&= \underset{(j+i)=odd}{\sum_{j=0}^{i}} \frac{(\frac{i+j-1}{2})!}{j! (\frac{i-j-1}{2})!} {}^C_0D^\alpha_x(x^j)\label{eq4.1}.
\end{align}
Now from the equation \eqref{eq1}, we have
\begin{align}
	{}^C_0D^\alpha_x(x^j) &= 0, \hspace{0.3cm} \text{If} \hspace{0.3cm} \lceil \alpha \rceil<j,\label{eq5.1}\\
	{}^C_0D^\alpha_x(x^j) & = \frac{j!}{\Gamma(j+1-\alpha)}x^{j-\alpha}, \hspace{0.3cm} \text{If} \hspace{0.3cm} \lceil \alpha \rceil \geq j\label{eq6.1}.
\end{align}
Now using above equations \eqref{eq5.1} and \eqref{eq6.1} in the equation \eqref{eq4.1}, we have
\begin{equation}
	\label{eq4}
	{}^C_0D^\alpha_x F_i(x) =
	\begin{cases}
		0\hspace{3.8cm}i={0,1,2,...,\lceil\alpha\rceil},\\
		\underset{(j+i)=odd}{\sum_{j=\lceil\alpha\rceil}^{i}} \frac{(\frac{i+j-1}{2})!}{(\frac{i-j-1}{2})!(j-\alpha)!} x^{j-\alpha},\hspace{0,2cm}  i\in N,\hspace{0,2cm}i\geq \lceil\alpha\rceil,
	\end{cases}
\end{equation}

\section{Architecture of FNN}
\begin{figure}[H]
	\centering
	\includegraphics[width=0.7\textwidth]{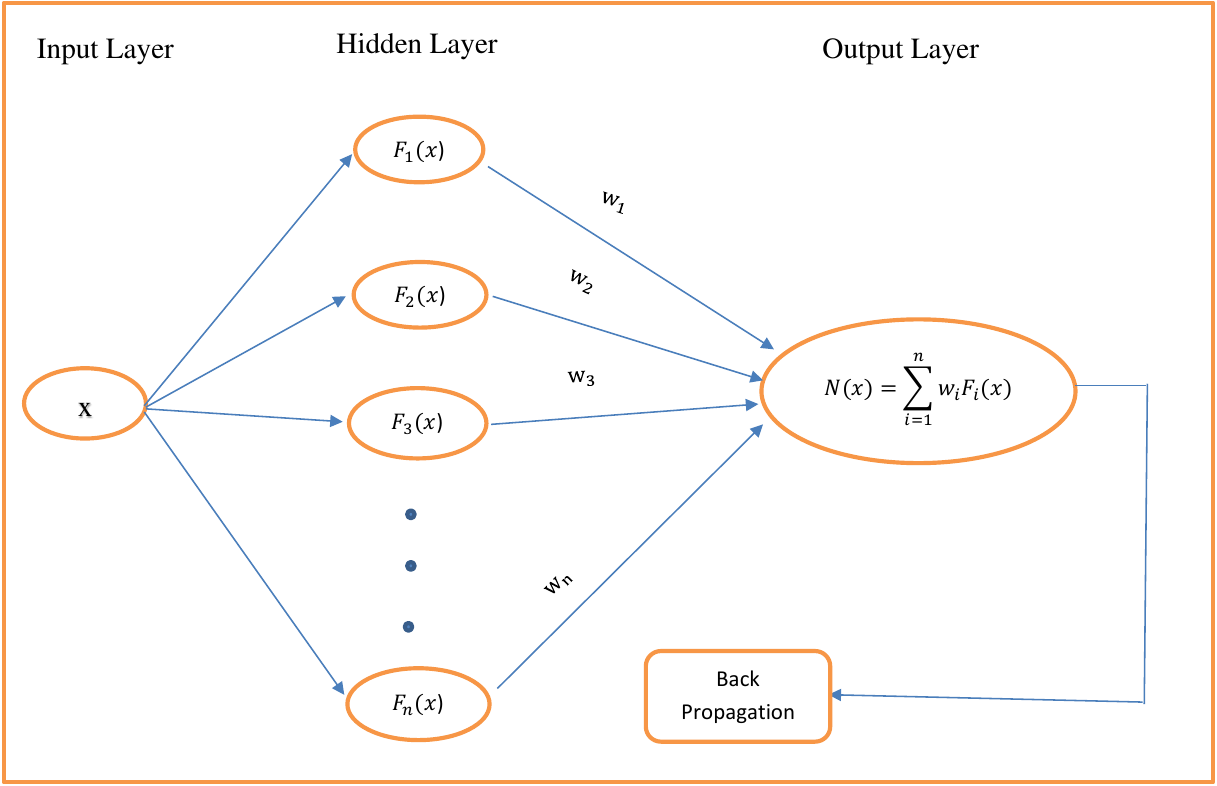}
	\caption{Structure of FNN to solve FDEs}
	\label{pic1}
\end{figure}

 Figure \ref{pic1} depicts the architecture of the FNN, which consists of input, hidden layer, and output layers. Input and output layers have one, and the hidden layer has $n$ numbers of perceptrons. The quantity of perceptrons in the hidden layer, denoted by $n$, varies based on the specific FDE being addressed for solution. Various degrees of Fibonacci polynomial are being used as an activation function in the hidden layer of the FNN. FNN will operate in the following ways:\\
\\
$\bullet$ From input layer, a input $x$ will be given.\\
$\bullet$ In the hidden layer the inputs from the input layer is activated through the different degrees of Fibonacci polynomials.\\
$\bullet$ The output layer is formed by combining inputs from the preceding layer using distinct weights $w$.

\section{Method to use FNN to solve FDEs}
In this section, the authors elaborate on a method for solving the general form of fractional order differential equations using the FNN. Let us consider the general form of FDE
\begin{equation}
	\label{eq5}
	{}^C_0D_x^{\alpha_m}y(x)=f(x,y(x),{}^C_0D_x^{\alpha_1}y(x),{}^C_0D_x^{\alpha_2}y(x),...,{}^C_0D_x^{\alpha_{m-1}}y(x)),\hspace{0.4cm}\alpha_m>\alpha_{m-1},...,\alpha_1>0
\end{equation}
with the ICs, $y^{(k)}(0)=y_0^k, \hspace{0.3cm}\text{For}\hspace{0.2cm}k=0,1,2...,\lceil\alpha\rceil$. The authors have considered the above problem with IC for the demonstration of the method but this technique can be used to solve the above problem with any type of conditions. Let's denote the trial solution of the FDE \eqref{eq5} as $u_{t}(t)=N(x)$, where $N(x)$ represents the output of the FNN and is defined as

\begin{equation}
	\label{eq6}
	N(x)=\sum_{i=i}^{n}w_iF_i(x).
\end{equation}
where, $w_i$'s are weights which are later be updated during the training of the FNN with appropriate method. Now the problem \eqref{eq5} will be transformed in to the following equation on using \eqref{eq6} in it.
\begin{equation}
	\label{eq7}
	{}^C_0D_x^{\alpha_m}N(x)=f(x,N(x),{}^C_0D_x^{\alpha_1}N(x),{}^C_0D_x^{\alpha_2}N(x)),
\end{equation}
with the following initial conditions
\begin{equation}
	\label{eq8}
	N^{(k)}(0)=y_0^k, \hspace{0.3cm}\text{For}\hspace{0.2cm}k=0,1,2...,\lceil\alpha\rceil,
\end{equation}
where,
\begin{equation*}
	{}^C_0D_x^{\alpha_m}N(x)=\sum_{i=1}^{n}w_i\Big({}^C_0D_x^{\alpha_m}F_i(x)\Big).
\end{equation*}
${}^C_0D_x^{\alpha_m}F_i(x)$ can easily be calculated from the equation  \eqref{eq4}. Now to train the FNN, the authors have transformed the above problem \eqref{eq7} with the initial conditions \eqref{eq8} in to the following minimization problem.
\begin{equation}
	\label{eq9}
	\underset{w}{min}\Bigl\{\sum_{x_p}\Big(	{}^C_0D_x^{\alpha_m}N(x)-f(x,N(x),{}^C_0D_x^{\alpha_1}N(x),{}^C_0D_x^{\alpha_2}N(x)) \Big)^2+\sum_{k=0}^{\lceil\alpha\rceil}\Big(N^{(k)}-y_0^k\Big)^2  \Bigr\}
\end{equation}
where, $x_p$'s the training points. Minimizing the above cost function \eqref{eq9} with respect to weights $w$'s is the training of neural network. In the forthcoming section, Marquardt's method has been discussed to update the values of weights $w$'s to minimize the cost function.

\section{Training of FNN}
Through the combination of  Newton's and Cauchy's optimization method, we can obtain Marquardt's \cite{ravindran2006engineering}. This is a second order method because it consist of Newton's method. The working procedure of Marquardt's method is as fallows \\

\begin{tikzpicture}	[node distance = 1cm,
	state/.style = {%
		rectangle,
		draw = black,
		inner sep = 0pt,
		minimum size = 10mm,
		thick,
	},
	beta/.style = {%
		node distance = 2mm,
		inner sep = 1pt,
	},
	auto,
	]
	
	\node (a1) [align=center,draw, terminal, minimum width=1cm, minimum height=1cm] {Assign $w^{(0)}$ \\
		(Choose the entries of $w^{(0)}$ non-identical and random)\\
		$M$ = Maximum number of iteration\\
		$\epsilon$ = convergence criteria\\
		Take $k=0,\hspace{0.2cm}\lambda^{(0)}=10^4$\\};
	
	\node (a2) [align=center, below=of a1, draw, rectangle, minimum width=3cm, minimum height=2cm] {
		Calculate $\nabla E(w^{(k)})$};
	
	\node (a3) [below= of a2, draw, rectangle, minimum width=3cm, minimum height=2cm] {If $E(w^{(k)})<\epsilon$? and $M\leq k$};
	
	\node (a4) [align=center, below= of a3, draw, rectangle, minimum width=3cm, minimum height=2cm] {Calculate \\
		$s(w^{(k)})=\big[H^{(k)}+\lambda^{(k)}I\big]^{-1}\nabla E(w^{(k)})$\\
		and $w^{(k+1)}=w^{(k)}-s(w^{(k)})$};

	\node (a5) [align=center, below= of a4, draw, rectangle, minimum width=3cm, minimum height=2cm] {If $E(w^{(k)})>E(w^{(k+1)})$?};
	
	\node (a6) at (-4,-16) [align=center, draw, rectangle, minimum width=3cm, minimum height=2cm] { Take \\
		$\lambda^{(k+1)}=\lambda^{(k)}/4$ and $k=k+1$};
	
	\node (a7) at (5,-16) [align=center, draw, rectangle, minimum width=3cm, minimum height=2cm] { $\lambda^{(k)}=2\lambda^{(k)}$};
	
	\node(a8) [align=center,right = of a3, draw, rectangle, minimum width=3cm, minimum height=2cm] {Print the result and stop};

	\draw [->] (a1) to node {} (a2);
	\draw [->] (a2) to node {} (a3);
	\draw [->] (a3) to node {Yes} (a4);
	\draw [->] (a3) to node {No} (a8);
	\draw [->] (a4) to node {Yes} (a5);
	\draw [->] (a5) to node {Yes} (a6);
	\draw [->] (a5) to node {No} (a7);
	\draw [->] (a6.north)..controls +(left:5mm) and +(left:50mm).. (a2.west);
	
	\draw [->] (a7.north)..controls +(right:5mm) and +(right:30mm).. (a4.east);
\end{tikzpicture}
\\
\\
Where, $E(w^{(k)})$, $\nabla E(w^{(k)})$ and $H^{(k)}$ denotes the values at $k^{th}$ iteration.
\begin{equation*}
	\begin{split}
		&E(w^{(k)})=\frac{1}{2n}\sum_{x_p}\Big(	{}^C_0D_x^{\alpha_m}N^{(k)}(x)-f(x,N^{(k)}(x),{}^C_0D_x^{\alpha_1}N^{(k)}(x),{}^C_0D_x^{\alpha_2}N^{(k)}(x)) \Big)^2+\frac{1}{2}\sum_{k=0}^{\lceil\alpha\rceil}\Big(N^{(k)}-y_0^{(k)}\Big)^2,
	\end{split}
\end{equation*}
\begin{small}
	\begin{equation*}
		\nabla E(w^{(k)})=\begin{bmatrix}
			\frac{\partial E(w^{(k)})}{\partial w_{1}}\\
			\frac{\partial E(w^{(k)})}{\partial w_{2}}\\
			\vdots\\
			\frac{\partial E(w^{(k)})}{\partial w_{l}}\\
			\vdots\\
			\frac{\partial E(w^{(k)})}{\partial w_{n}}
		\end{bmatrix},\hspace{0.2cm}H^{(k)}=\begin{bmatrix}
			\frac{\partial^2 E(w^{(k)})}{\partial w_{1}^2}&\frac{\partial^2 E(w^{(k)})}{\partial w_{1}\partial w_{2}}& \cdot&\frac{\partial^2 E(w^{(k)})}{\partial w_{1}\partial w_{n}}\\
			\vdots&\vdots&\cdot&\vdots\\
			\frac{\partial^2 E(w^{(k)})}{\partial w_{l}\partial w_{1}}&\frac{\partial^2 E(w^{(k)})}{\partial w_{l}\partial w_{2}}& \cdot&\frac{\partial^2 E(w^{(k)})}{\partial w_{l}\partial w_{n}}\\
			\vdots&\vdots&\cdot&\vdots\\
			\frac{\partial^2 E(w^{(k)})}{\partial w_{n}\partial w_{1}}&\frac{\partial^2 E(w^{(k)})}{\partial w_{n}\partial w_{2}}& \cdot&\frac{\partial^2 E(w^{(k)})}{\partial w_{n}^2}
		\end{bmatrix},
	\end{equation*}	
\end{small}
where,
\begin{align}
	\frac{\partial E(w^{(k)})}{\partial w_{l}} & =\frac{1}{n} \sum_{x_p}\Big(	{}^C_0D_x^{\alpha_m}N^{(k)}(x)-f(x,N^{(k)}(x),{}^C_0D_x^{\alpha_1}N^{(k)}(x),{}^C_0D_x^{\alpha_2}N^{(k)}(x)) \Big) \Big(\frac{\partial \big({}^C_0D_x^{\alpha_m}N^{(k)}(x)\big)}{\partial w_l}\nonumber\\
	&-\frac{\partial f(x,N^{(k)}(x),{}^C_0D_x^{\alpha_1}N^{(k)}(x),{}^C_0D_x^{\alpha_2}N^{(k)}(x))}{\partial w_l} \Big)+\sum_{k=0}^{\lceil\alpha\rceil}\Big(N^{(k)}-y_0^k\Big)\frac{\partial N^{(k)}}{\partial w_l},\nonumber\\
	\frac{\partial^2 E(w^{(k)})}{\partial w_{l}w_{m}} & = \frac{1}{n} \sum_{x_p}\Big(\frac{\partial \big({}^C_0D_x^{\alpha_m}N^{(k)}(x)\big)}{\partial w_m}
	-\frac{\partial f(x,N^{(k)}(x),{}^C_0D_x^{\alpha_1}N^{(k)}(x),{}^C_0D_x^{\alpha_2}N^{(k)}(x))}{\partial w_m} \Big) \Big(\frac{\partial \big({}^C_0D_x^{\alpha_m}N^{(k)}(x)\big)}{\partial w_l}\nonumber\\
	&-\frac{\partial f(x,N^{(k)}(x),{}^C_0D_x^{\alpha_1}N^{(k)}(x),{}^C_0D_x^{\alpha_2}N^{(k)}(x))}{\partial w_l} \Big)+\sum_{k=0}^{\lceil\alpha\rceil}\frac{\partial N^{(k)}}{\partial w_m}\frac{\partial N^{(k)}}{\partial w_l}.\nonumber
\end{align}

In the subsequent section, the authors employed the discussed method on five examples having exact solutions to verify its accuracy numerically. They further demonstrated its superiority over previously established methods by comparing the absolute errors obtained with those of the proposed approach and the earlier methods.

\section{Numerical Examples}
 
In this section, the authors showed the reliability of the discussed numerical technique by using it to solve FDES having the exact solutions. Further, compared the numerical results obtained from this method and the previously solved method to reveal the reliability of the discussed method. The authors use the Python-3.7.9 version to use the discussed method in the following examples. \par 

For Python code of following example visit "\url{https://github.com/Kushaldhardwivedi/Differential-Equation-articles}" \\
\\
\textbf{Example 1.} The subsequent FDE
\begin{equation}
	{}_0^CD^{\alpha}_ty(t)+y(t)=f(t),\hspace{1cm}0<\alpha\leq 1,
\end{equation}
with IC $y(0)=0$ have the exact solution $y(t)=t^2$, where
\begin{equation*}
	f(t)=t^2+\frac{2t^{2-\alpha}}{\Gamma(3-\alpha)}.
\end{equation*}

this problem  has been solved with $n=3$ perceptrons in the hidden layer of the FNN and trained the FNN with $10$ training points of the domain $0\leq x\leq1$ over $18$ iterations. Further, compared the numerical results obtained by solving this problem with the suggested technique in the following Table \ref{table:1}.
\begin{table}[H]
	\begin{center}
		\caption{Comparison of absolute error of example 1 with different method.}
		\label{table:1}
		\begin{tabular}{c c c c c c c c c }
			\hline
			&$\alpha=0.25$& & &$\alpha=0.5$& & &$\alpha=0.75$ &\\
			\cline{2-3} \cline{5-6} \cline{8-9}
			&Proposed 	&Method & &Proposed& Method & &Proposed& Method \\
			$t$&method &\cite{shah2017numerical} & & Method & \cite{shah2017numerical} & & Method & \cite{shah2017numerical}\\
			\hline
			0.1	&6.93889e-18 &0.0000e-4 & &6.93889e-18 &0.0000e-4 & &1.04083e-16 &0.0000e-4 \\
			0.2	&2.77556e-17 &1.0000e-4 & &2.77556e-17 &1.0000e-4 & &2.77556e-17
			&0.0000e-4 \\
			0.3	&5.55112e-17 &1.0000e-4 & &5.55112e-17 &1.0000e-4 & &1.66533e-16
			&0.0000e-4 \\
			0.4	&1.11022e-16 &0.0000e-4 & &1.11022e-16 &0.0000e-4 & &2.22045e-16
			&0.0000e-4 \\
			0.5	&0.00000        &3.0000e-4 & &2.22045e-16 &3.0000e-4 & &0.0000 &2.0000e-4 \\
			0.6	&0.000000        &1.0000e-4 & &2.22045e-16 &1.0000e-4 & &1.11022e-16
			&0.0000e-4 \\
			0.7	&1.11022e-16 &1.0000e-4 & &1.11022e-16 &0.0000e-4 & &1.11022e-16
			&1.0000e-4 \\
			0.8	&0.00000        &1.0000e-4 & &0.00000        &1.0000e-4 & &1.11022e-16
			&0.0000e-4 \\
			0.9	&0.00000        &0.0000e-4 & &0.00000        &0.0000e-4 & &0.00000 &0.0000e-4 \\
			\hline
		\end{tabular}	
	\end{center}	
\end{table}

\begin{figure}[H]
	\centering
	\includegraphics[width=0.7\textwidth]{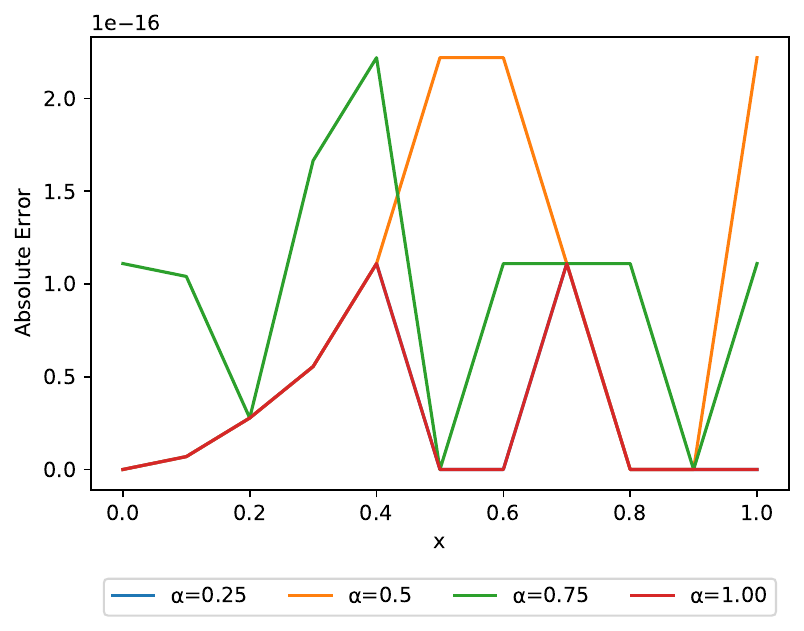}
	\caption{Absolute Error for distinct values of $\alpha$ of Example 1.}
	\label{pic:1}
\end{figure}
This problem has also been solved by many other authors for different values of $\alpha$, viz., Kai Diethelm \cite{diethelm1997algorithm} solved this problem with the best absolute error of order $e-4$, Shah et al. \cite{shah2017numerical} solve this problem by employing HaarWavelet Operational Matrix Method with the lowest error of order $e-4$ and Li \cite{li2012numerical} solved this problem with the minimum error of order $e-12$, etc. From Table \ref{table:1} and Fig \ref{pic:1}, it is clearly visible that our method is performing far better than previously solved methods.\\
\\
\textbf{Example 2.} The following differential equation
\begin{equation}
	\frac{d^2y(t)}{dt^2}+{}^C_0D^{\frac{1}{2}}_ty(t)+y=t^3+6t+\frac{3.2}{\Gamma(0.5)}t^{2.5}
\end{equation}

with ICs $y(0)=0$, $y'(0)=0$ have the exact solution $y(t)=t^3$. Authors have solved this problem with $n=4$ perceptrons in the hidden layer of the FNN and trained the FNN with $10$ training points of the domain $0\leq t\leq 1$  over $17$ iterations.

\begin{table}[H]
	\begin{center}
		\caption{Maximum absolute error of example 2.}
		\label{table:2}
		\begin{tabular}{c c c c }
			\hline
			& Numerical & Exact & Error\\
			$t$	& Solution & Solution &  \\
			\hline
			0.0 &3.94594e-17& 0.00     & 3.94594e-17\\
			0.1 &  0.001     &  0.001   & 5.48606e-17\\
			0.2 &  0.008     &0.008     & 5.0307e-17      \\
			0.3 &  0.027     &0.027     & 1.04083e-17     \\
			0.4 &  0.064     & 0.064    & 4.16334e-17   \\
			0.5 &  0.125     & 0.125    & 0.00\\
			0.6 &  0.216     & 0.216    & 8.32667e-17 \\
			0.7 &  0.343     & 0.343    & 1.11022e-16  \\
			0.8 &  0.512     & 0.512    & 1.11022e-16   \\
			0.9 &  0.729     & 0.729    & 2.22045e-16    \\
			1.0 &  1.000     & 1.000    & 0.00     \\
			\hline
		\end{tabular}	
	\end{center}	
\end{table}
\begin{figure}[H]
	\centering
	\includegraphics[width=0.7\textwidth]{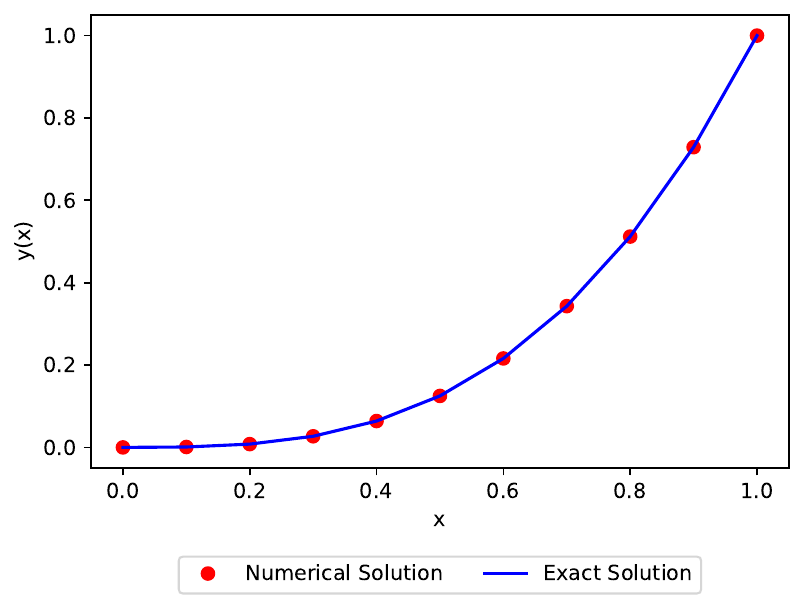}
	\caption{Comparison of numerical and exact solution of Example 2.}
	\label{pic:2}
\end{figure}
\begin{figure}[H]
	\centering
	\includegraphics[width=0.7\textwidth]{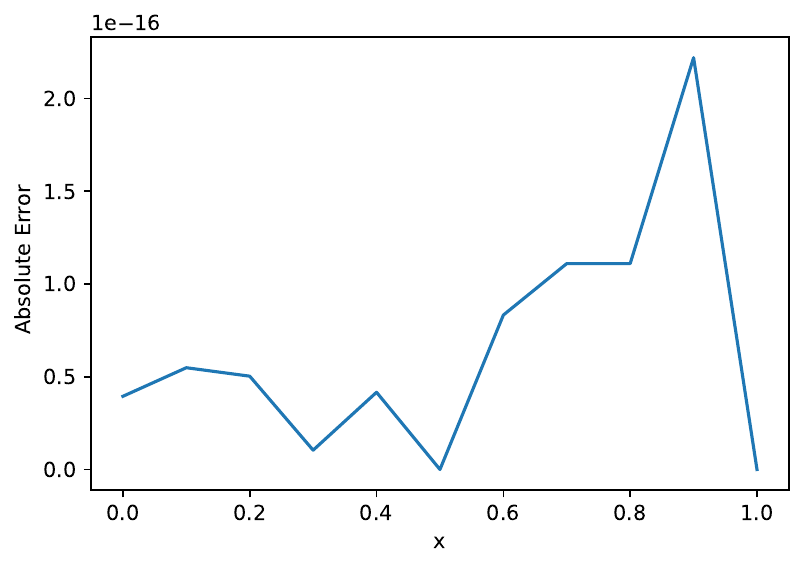}
	\caption{Absolute Error of Example 2.}
	\label{pic:3}
\end{figure}
In Table 2, the authors compared the numerical solution obtained from the discussed and exact solution of the problem. It is clearly visible  From Table \ref{table:2} that the proposed method is performing accurately. Ford et al. \cite{ford2009systems} solved this problem with Systems-based decomposition schemes with the best error of order $e-7$. Shiralashetti et al. \cite{shiralashetti2016efficient} also solved this problem using Haar wavelet collocation method with best error of order $e-9$. From Fig \ref{pic:3}, it is clearly visible that the proposed method is performing far better than previous methods.\\
\\
\textbf{Example 3.} The following nonlinear differential equation
\begin{equation}
	{}^C_0D^{2.2}_ty(t)+{}^C_0D^{0.75}_ty(t)+{}^C_0D^{1.25}_ty(t)+y^3=f(t),
\end{equation}
with the ICs $y(0)=0$, $y'(0)=0$ and $y''(0)=0$ have the exact solution $y(t)=\frac{t^3}{3}$, where
\begin{equation*}
	f(t)=\frac{2t^{0.8}}{\Gamma(1.8)}+\frac{2t^{2.25}}{\Gamma(3.25)}+\frac{2t^{1.75}}{\Gamma(2.75)}+\frac{t^9}{27}.
\end{equation*}
The authors have solved this problem with $n=4$ perceptrons in the hidden layer of FNN and trained the FNN with $10$ training point of the domain $0\leq t\leq1$ over $16$ iterations.
\begin{table}[H]
	\begin{center}
		\caption{Comparison absolute error of example 2 with distinct methods.}
		\label{table:3}
		\begin{tabular}{c c c c }
			\hline
			& Proposed & \cite{li2010haar} & \cite{shiralashetti2016efficient}\\
			$t$	& Method & Method & Method  \\
			\hline
			0.1 &1.87025e-17&1.530e-05  &5.052e-06\\
			0.2 &1.0842e-17 &3.690e-05  &1.020e-05 \\
			0.3 &3.1225e-17 &5.420e-05  &1.478e-05 \\
			0.4 &3.1225e-17 &5.91e-05  &1.871e-05 \\
			0.5 &2.08167e-17&9.10e-05  &2.197e-05 \\
			0.6 &8.32667e-17&8.28e-05  &2.459e-05 \\
			0.7 &4.16334e-17&1.140e-04  &2.660e-05 \\
			0.8 &8.32667e-17&1.263e-04  &2.808e-05 \\
			0.9 &2.77556e-17&1.119e-04  &2.914e-05 \\
			\hline
		\end{tabular}	
	\end{center}	
\end{table}
\begin{figure}[H]
	\centering
	\includegraphics[width=0.7\textwidth]{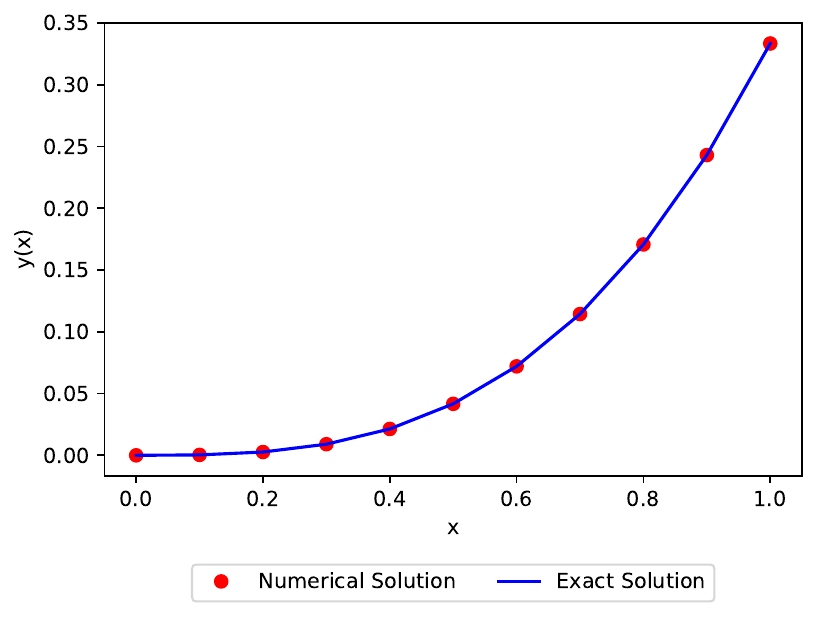}
	\caption{Comparison of numerical and exact solution of Example 3.}
	\label{pic:4}
\end{figure}
\begin{figure}[H]
	\centering
	\includegraphics[width=0.7\textwidth]{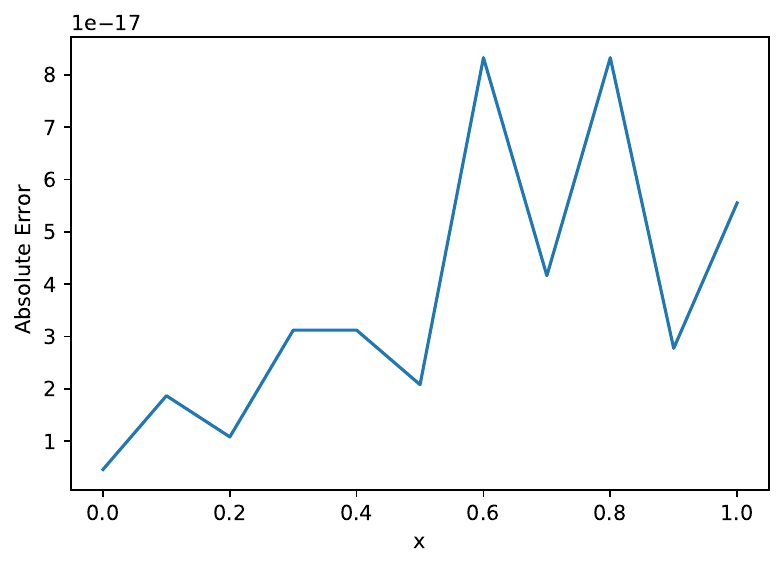}
	\caption{Absolute Error of Example 3.}
	\label{pic:5}
\end{figure}
From Fig \ref{pic:4} and Fig \ref{pic:5} it clear that the discussed method is giving very accurate result. Many authors have also solved this problem with the different methods viz., El-Mesiry et al. \cite{el2005numerical} solve this problem with the best error of order $e-3$, Li et al. \cite{li2010haar}, and Shiralashetti et al. \cite{shiralashetti2016efficient} have also solved this problem with Haar wavelet. This comparison between numerical results obtained from these methods and our from our discussed method is depicted through Table \ref{table:3}. From Table \ref{table:3} it is clearly visible that our method is performing far more accurately than these methods.\\
\\
\textbf{Example 4.} The following FDE
\begin{equation}
	{}^C_0D^{\alpha}_ty(t)+y(t)=f(t),
\end{equation}
with IC $y(0)=1$, have the exact solution $y(t)=1-4t+5t^2$, where
\begin{equation*}
	f(t)=1-4t+5t^2-\frac{4}{\Gamma(2-\alpha)}t^{1-\alpha}+\frac{10}{\Gamma(3-\alpha)}t^{2-\alpha}.
\end{equation*}
\begin{table}[H]
	\begin{center}
		\caption{Maximum absolute error of example 4.}
		\label{table:4}
		\begin{tabular}{c c c c c c c c c }
			\hline
			&$\alpha=0.25$& & &$\alpha=0.5$& & &$\alpha=0.75$ &\\
			\cline{2-3} \cline{5-6} \cline{8-9}
			&Proposed 	&Method & &Proposed& Method & &Proposed& Method \\
			$t$&method &\cite{shah2017numerical} & & Method & \cite{shah2017numerical} & & Method & \cite{shah2017numerical}\\
			\hline
			0.1	&1.22124e-15 &6.0000e-4 & &5.55111e-16 &1.0000e-4 & &3.33066e-16
			&6.0000e-4 \\
			0.2	&8.32667e-15 &4.0000e-4 & &5.55111e-16 &1.0000e-4 & &2.10942e-15
			&4.0000e-4 \\
			0.3	&1.33781e-14 &2.0000e-4 & &5.55111e-17 &1.0000e-4 & &3.60822e-15
			&2.0000e-4 \\
			0.4	&1.69864e-14 &2.0000e-4 & &9.99200e-16 &2.0000e-4 & &4.55191e-15
			&2.0000e-4 \\
			0.5	&1.59872e-14 &5.0000e-4 & &0.00000e+00 &9.0000e-4 & &4.44089e-15
			&5.0000e-4 \\
			0.6	&1.26565e-14 &3.0000e-4 & &6.66133e-16 &5.0000e-4 & &3.77475e-15
			&3.0000e-4 \\
			0.7	&8.88178e-15 &1.0000e-4 & &0.00000e+00 &0.0000e-4 & &3.55271e-15
			&1.0000e-4 \\
			0.8	&1.33226e-15 &1.0000e-4 & &4.44089e-16 &1.0000e-4 & &2.22044e-15
			&1.0000e-4 \\
			0.9	&9.32587e-15 &1.0000e-4 & &1.33226e-15 &0.0000e-4 & &1.33226e-15
			&1.0000e-4 \\
			\hline
		\end{tabular}	
	\end{center}	
\end{table}
\begin{figure}[H]
	\centering
	\includegraphics[width=0.7\textwidth]{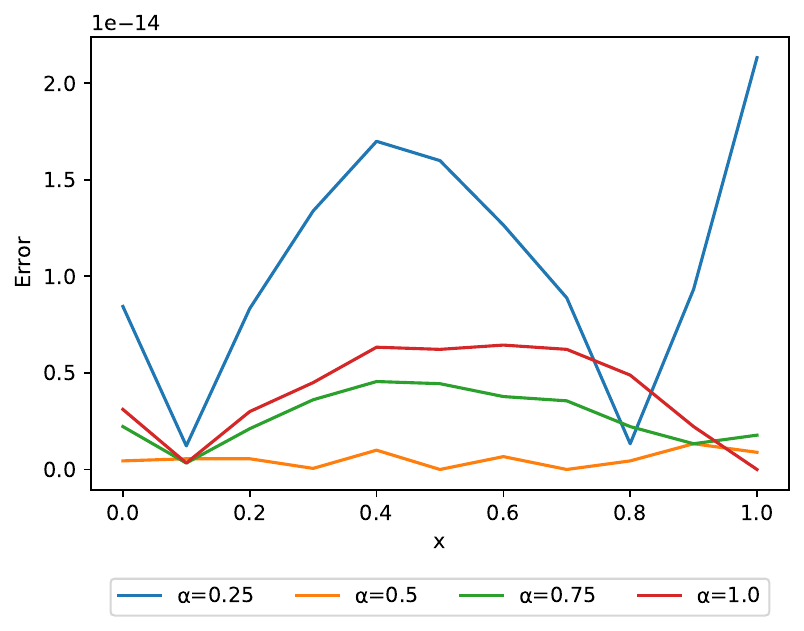}
	\caption{Absolute Error of Example 4.}
	\label{pic:6}
\end{figure}
The authors have solved this problem with $n=3$ perceptrons in the hidden layer of FNN and trained the neural network on $10$ training points of the domain $0\leq x\leq 1$ over $28$ iterations. It can be observed from Table \ref{table:4} that the discussed method gives far more reliable results than the compared method.\\
\\
\textbf{Example 5.} Let us consider the following initial value problem of Bagley–
Torvik equation

\begin{table}[H]
	\begin{center}
		\caption{Comparison between absolute errors of example 5 with different methods.}
		\label{table:5}
		\begin{tabular}{c c c c }
			\hline
			& Method &Method & Proposed\\
			$t$	&\cite{raja2017design} & \cite{verma2020numerical} & Method \\
			\hline
			0.0 & 7.63e-06    & 0.00     & 0.00\\
			0.1 &  8.73e-06   & 7.59e-08   & 6.93889e-18\\
			0.2 &  1.12e-06   &1.57e-07     & 2.77556e-17      \\
			0.3 &  1.08e-06   &2.09e-07     & 5.55112e-17 \\
			0.4 &  8.19e-06   & 2.74e-07   & 1.11022e-16\\
			0.5 &  7.06e-06   & 3.61e-07   &0.00\\
			0.6 &  1.01e-05   & 4.37e-07    & 0.00 \\
			0.7 &  1.60e-05   & 4.68e-07    &1.11022e-16 \\
			0.8 &  2.03e-05   & 4.64e-07   &   0.00\\
			0.9 &  1.86e-05   & 4.87e-07    &  0.00   \\
			1.0 &  1.24e-05   & 5.74e-07    &  0.00    \\
			\hline
		\end{tabular}	
	\end{center}	
\end{table}
\begin{equation}
	\frac{d^2y(t)}{dt^2}+{}^C_0D^{\alpha}_ty(t)+y(t)=2+4\sqrt{\frac{t}{\pi}}+t^2,\hspace{0.2cm}0\leq t\leq1,
\end{equation}
with the initial conditions $y(0)=0$ and $y'(0)=0$. This problem have the exact solution $y(t)=t^2$ for $\alpha=1.5$. This problem has been solved by taking $n=4$ perceptron in the hidden layers of the FNN with $10$ training points of the domain $0\leq t\leq 1$ over $16$ iterations.

\begin{figure}[H]
	\centering
	\includegraphics[width=0.7\textwidth]{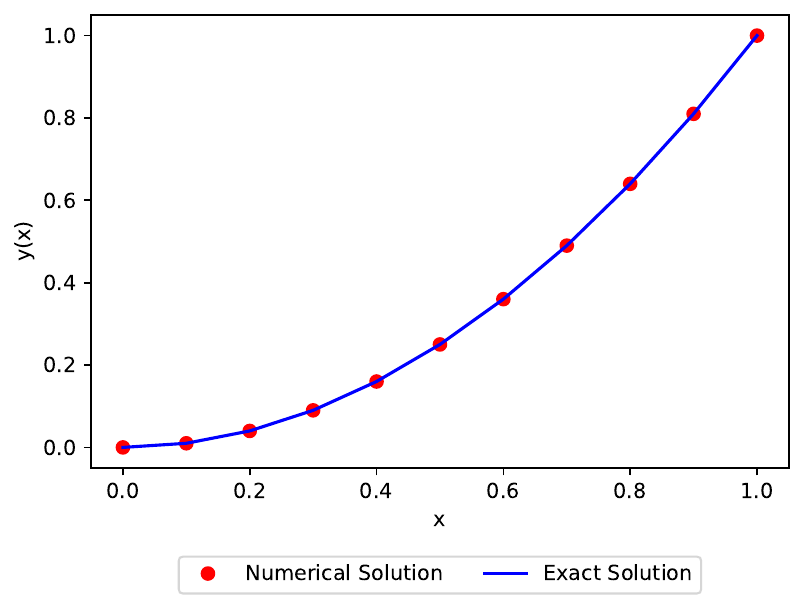}
	\caption{Comparison of numerical and exact solution of Example 5 for $\alpha=1.5$.}
	\label{pic:7}
\end{figure}
\begin{figure}[H]
	\centering
	\includegraphics[width=0.7\textwidth]{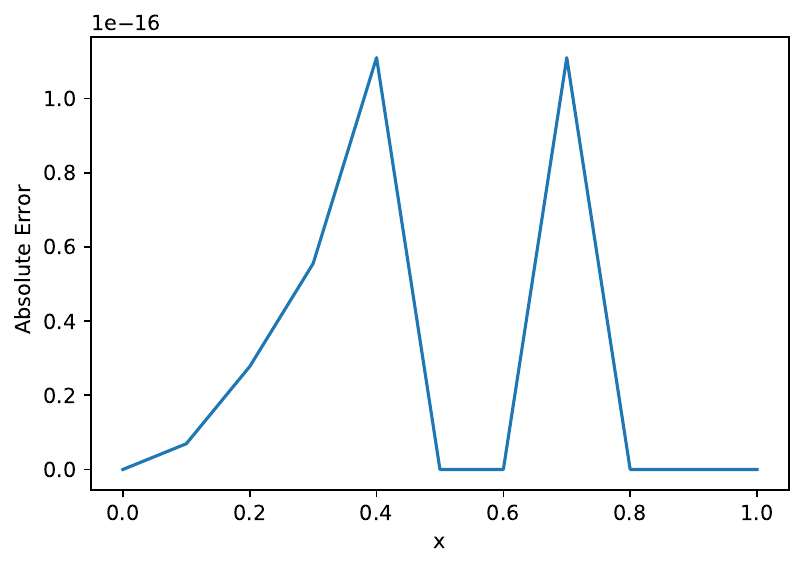}
	\caption{Absolute Error of Example 5 for $\alpha=1.5$.}
	\label{pic:8}
\end{figure}
From the table \ref{table:5} it can be clearly observe that the proposed method is performing far better than the method discussed in \cite{raja2017design} and \cite{verma2020numerical}. Figures \ref{pic:7} and \ref{pic:8} have been plotted to compare the numerical and exact solution and to observe the absolute error respectively after using the discussed method on Example 5 for $\alpha=1.5$.

\section{Conclusion}
In this paper, the authors used three layers of the FNN to solve the FDEs. In the second layer, different orders of the Fibonacci polynomial has been used as an activation function. To solve the desired FDEs, the authors used the FNN as a trial solution of the problem and then trained the weights by Marquardt's method. After discussing the method, the authors used it in four examples to reveal the efficiency of the proposed method. Furthermore, the author showed that the discussed approach is giving more reliable numerical outcomes than the previously solving method.

\section{Conflict of interest}
The authors have no conflict of interest to mention.

\bibliographystyle{elsarticle-num} 
\bibliography{reference}

\begin{thebibliography}{10}
\expandafter\ifx\csname url\endcsname\relax
  \def\url#1{\texttt{#1}}\fi
\expandafter\ifx\csname urlprefix\endcsname\relax\def\urlprefix{URL }\fi
\expandafter\ifx\csname href\endcsname\relax
  \def\href#1#2{#2} \def\path#1{#1}\fi

\bibitem{zhao2006new}
W.~Zhao, L.~Chen, S.~Peng, A new kind of accurate numerical method for backward
  stochastic differential equations, SIAM Journal on Scientific Computing
  28~(4) (2006) 1563--1581.

\bibitem{ur2012numerical}
M.~ur~Rehman, R.~A. Khan, A numerical method for solving boundary value
  problems for fractional differential equations, Applied Mathematical
  Modelling 36~(3) (2012) 894--907.

\bibitem{yan2014higher}
Y.~Yan, K.~Pal, N.~J. Ford, Higher order numerical methods for solving
  fractional differential equations, BIT Numerical Mathematics 54~(2) (2014)
  555--584.

\bibitem{li2017new}
X.~Li, H.~Li, B.~Wu, A new numerical method for variable order fractional
  functional differential equations, Applied Mathematics Letters 68 (2017)
  80--86.

\bibitem{baleanu2018chebyshev}
D.~Baleanu, B.~Shiri, H.~Srivastava, M.~Al~Qurashi, A chebyshev spectral method
  based on operational matrix for fractional differential equations involving
  non-singular mittag-leffler kernel, Advances in Difference Equations 2018~(1)
  (2018) 1--23.

\bibitem{shiri2019system}
B.~Shiri, D.~Baleanu, System of fractional differential algebraic equations
  with applications, Chaos, Solitons \& Fractals 120 (2019) 203--212.

\bibitem{shiri2020collocation}
B.~Shiri, G.-C. Wu, D.~Baleanu, Collocation methods for terminal value problems
  of tempered fractional differential equations, Applied Numerical Mathematics
  156 (2020) 385--395.

\bibitem{dadkhah2020spline}
E.~Dadkhah~Khiabani, H.~Ghaffarzadeh, B.~Shiri, J.~Katebi, Spline collocation
  methods for seismic analysis of multiple degree of freedom systems with
  visco-elastic dampers using fractional models, Journal of Vibration and
  Control 26~(17-18) (2020) 1445--1462.

\bibitem{dadkhah2020visco}
E.~Dadkhah, B.~Shiri, H.~Ghaffarzadeh, D.~Baleanu, Visco-elastic dampers in
  structural buildings and numerical solution with spline collocation methods,
  Journal of Applied Mathematics and Computing 63~(1) (2020) 29--57.

\bibitem{balasubramaniam2006solution}
P.~Balasubramaniam, J.~A. Samath, N.~Kumaresan, A.~V.~A. Kumar, Solution of
  matrix riccati differential equation for the linear quadratic singular system
  using neural networks, Applied Mathematics and Computation 182~(2) (2006)
  1832--1839.

\bibitem{reynaldi2012backpropagation}
A.~Reynaldi, S.~Lukas, H.~Margaretha, Backpropagation and levenberg-marquardt
  algorithm for training finite element neural network, in: 2012 Sixth
  UKSim/AMSS European Symposium on Computer Modeling and Simulation, IEEE,
  2012, pp. 89--94.

\bibitem{jafarian2017artificial}
A.~Jafarian, M.~Mokhtarpour, D.~Baleanu, Artificial neural network approach for
  a class of fractional ordinary differential equation, Neural Computing and
  Applications 28~(4) (2017) 765--773.

\bibitem{zuniga2017solving}
C.~Z{\'u}{\~n}iga-Aguilar, H.~Romero-Ugalde, J.~G{\'o}mez-Aguilar,
  R.~Escobar-Jim{\'e}nez, M.~Valtierra-Rodr{\'\i}guez, Solving fractional
  differential equations of variable-order involving operators with
  mittag-leffler kernel using artificial neural networks, Chaos, Solitons \&
  Fractals 103 (2017) 382--403.

\bibitem{podlubny1998fractional}
I.~Podlubny, Fractional differential equations: an introduction to fractional
  derivatives, fractional differential equations, to methods of their solution
  and some of their applications, Elsevier, 1998.

\bibitem{ravindran2006engineering}
A.~Ravindran, G.~V. Reklaitis, K.~M. Ragsdell, Engineering optimization:
  methods and applications, John Wiley \& Sons, 2006.

\bibitem{shah2017numerical}
F.~A. Shah, R.~Abass, L.~Debnath, Numerical solution of fractional differential
  equations using haar wavelet operational matrix method, International Journal
  of Applied and Computational Mathematics 3~(3) (2017) 2423--2445.

\bibitem{diethelm1997algorithm}
K.~Diethelm, An algorithm for the numerical solution of differential equations
  of fractional order, Electron. Trans. Numer. Anal 5~(1) (1997) 1--6.

\bibitem{li2012numerical}
X.~Li, Numerical solution of fractional differential equations using cubic
  b-spline wavelet collocation method, Communications in Nonlinear Science and
  Numerical Simulation 17~(10) (2012) 3934--3946.

\bibitem{ford2009systems}
N.~J. Ford, J.~A. Connolly, Systems-based decomposition schemes for the
  approximate solution of multi-term fractional differential equations, Journal
  of Computational and Applied Mathematics 229~(2) (2009) 382--391.

\bibitem{shiralashetti2016efficient}
S.~Shiralashetti, A.~Deshi, An efficient haar wavelet collocation method for
  the numerical solution of multi-term fractional differential equations,
  Nonlinear dynamics 83~(1-2) (2016) 293--303.

\bibitem{li2010haar}
Y.~Li, W.~Zhao, Haar wavelet operational matrix of fractional order integration
  and its applications in solving the fractional order differential equations,
  Applied Mathematics and Computation 216~(8) (2010) 2276--2285.

\bibitem{el2005numerical}
A.~El-Mesiry, A.~El-Sayed, H.~El-Saka, Numerical methods for multi-term
  fractional (arbitrary) orders differential equations, Applied Mathematics and
  Computation 160~(3) (2005) 683--699.

\bibitem{raja2017design}
M.~A.~Z. Raja, R.~Samar, M.~A. Manzar, S.~M. Shah, Design of unsupervised
  fractional neural network model optimized with interior point algorithm for
  solving bagley--torvik equation, Mathematics and Computers in Simulation 132
  (2017) 139--158.

\bibitem{verma2020numerical}
A.~Verma, M.~Kumar, Numerical solution of bagley--torvik equations using
  legendre artificial neural network method, Evolutionary Intelligence (2020)
  1--11.

\end{thebibliography}

\end{document}